\documentclass[a4paper]{article}
\usepackage{amsmath,amssymb,amsthm,amsfonts,amscd}
\theoremstyle{main}
\newtheorem{theorem}{Theorem}

\begin{document}
\author{Evgeny Shchepin}

 \title{ Riesz means and Greedy Sums.}
\maketitle

\paragraph{Motivation.}
The concept of greedy sum $\sum\limits_{s\in S} a_s$ of a numerical
array $\{a_s\}_{s\in S}$  was introduced in \cite{Shch}. 

The direct product $\{a_sb_t\}_{s,t\in S\times T}$ of two greedy
summable numeric arrays $\{a_s\}_{s\in S}$ $\{b_t\}_{t\in T}$ is
not greedy summable in general. 
To sum up such products we adopt the general method of Riesz means.

\paragraph{Arithmetic Riesz means.}
Let $\{\lambda_n\}$ be  an increasing to infinity sequence of
positive numbers and $\kappa$ be a positive number. The series
$\sum a_k$ is called \emph{summable} by Riesz method with
parameters $(\{\lambda_n\},\kappa$  to the sum $A$ if
\begin{equation}\label{r-l}
  A= \lim \limits_{\omega\to\infty} \sum\limits_{\lambda_n<\omega}
    \left(1-\frac{\lambda_n}\omega \right)^\kappa a_n,
\end{equation}

If we consider a sequence of complex numbers $a_n$, such that
$|a_n|\searrow 0$, then
  we can substitute in this definition
$\lambda_n=\frac1{|a_n|}$ and $\omega=\frac1\varepsilon$.
Riesz means with such parameters being applied to Dirichlet series
associated with a numerical array $\{a_s\}_{s\in S}$ (see \cite{Shch})
led us to the following definition.
 А numeric array
 $\{a_s\}_{s\in S}$ is called \emph{greedy $\kappa$-summable} for
 some positive  $\kappa$, if there exists the following limit \begin{equation}
\label{kappa-sum}\lim\limits_{\varepsilon\to0}
\sum\limits_{|a_s|>\varepsilon}a_s\left(1-\frac\varepsilon{|a_s|}\right)^\kappa,
\end{equation}
which we name as \emph{greedy $\kappa$-sum} of the array. For
$\kappa=0$ the concept of greedy $\kappa$-sum converts into the
original concept of greedy sum.

\begin{theorem}\label{riesz-compatibility}
If an array $\{a_s\}_{s\in S}$ has greedy $\kappa$-sum equal to
$A$ for some $\kappa$ then it has the greedy $\kappa'$-sum for any
$\kappa'<\kappa$.
\end{theorem}

The above theorem follows from the theorem 16 of \cite{H-R}.

\begin{theorem}\label{riesz-prod}
If an array $\{a_s\}_{s\in S}$ has greedy $\kappa_1$-sum equal to
$A$ and another array $\{b_t\}_{t\in T}$ has greedy $\kappa_2$-sum
equal to $B$, then their product $\{a_sb_t\}_{s,t\in S\times T}$
has greedy $(\kappa_1+\kappa_2+1)$-sum equal to $AB$.
\end{theorem}

 The last theorem follows  from  the theorem 56 of \cite{H-R}.

\end{document}